\documentclass[11 pt]{amsart}

\usepackage{amscd, amssymb}

\newtheorem{pr}{Proposition}
\newtheorem{lm}{Lemma}

\newtheorem{conj}{Conjecture}

\newcommand{\proj}{\mathbf P}

\newcommand{\rarr}{\rightarrow}
\newcommand{\oh}{{\mathcal{O}}}
\newcommand{\com}{\mathbb{C}}
\newcommand{\Q}{\mathbb{Q}}

\newcommand{\lan}{\langle}
\newcommand{\ran}{\rangle}
\newcommand{\eqq}{\cong}

\newcommand{\br}{{br}}

\def\scup{\mathbin{\text{\scriptsize$\cup$}}}

\newcommand{\bpf}{\noindent {\em Proof.} }
\newcommand{\epf}{\qed \vspace{+10pt}}

\newcommand{\hodge}{\mathbb{E}}

\newcommand{\ch}{\text{ch}}

\newcommand{\UM}{\mathcal{C}}
\newcommand{\sym}{{\text{Sym}}}

\newcommand{\mbar}[1]{{\overline{M}}_{#1}}
\newcommand{\lmb}[1]{\lambda_{#1}}
\newcommand{\kp}[1]{\kappa_{#1}}
\newcommand{\ta}[1]{\tau_{#1}}
\newcommand{\pa}{\partial}

\begin{document}
\title{Logarithmic series and 
Hodge integrals in the tautological ring}
\author{C. Faber and R. Pandharipande (with an appendix by D. Zagier)}
\date{9 March 2000}
\dedicatory{Dedicated to William Fulton on the occasion of his 60th
birthday}
\maketitle

\pagestyle{plain}
\setcounter{section}{-1}
\section{\bf{Introduction}}
\subsection{Overview}
Let $X_g$ be a nonsingular curve of genus $g\geq 2$ over
$\com$.
$X_g$ determines a point 
$[X_g] \in \overline{M}_g$ in the moduli space
of Deligne-Mumford stable genus $g$ curves.
The study of the Chow ring of the moduli space 
of curves was initiated by D. Mumford in [Mu].
In the past two decades, many remarkable properties
of these intersection rings have been discovered.
Our first goal in this paper is to describe 
a new perspective on the intersection theory
of the moduli space of curves
which encompasses advances
from both classical degeneracy studies and
topological gravity.
This approach is developed in Sections \ref{vvv} - \ref{wwww}
of the Introduction.

The main new results of the paper are computations of basic Hodge integral
series in $A^*(M_g)$ encoding the canonical evaluations of
$\kappa_{g-2-i} \lambda_i$. The motivation for the study of these
tautological elements and the series results are given in  
Section \ref{sss} of the Introduction. 
The body of the
paper contains the Hodge integral derivations.

\subsection{Moduli filtration}
\label{vvv}
Let $X_g$ be a fixed nonsingular curve.
We will consider the moduli filtration:
\begin{equation}
\label{fff}
\overline{M}_g \supset M^c_g \supset  M_g \supset \{ [X_g] \}.
\end{equation}
Here, $M_g$ is the moduli space of nonsingular genus $g$
curves, and $M_g^c$ is the moduli space of 
stable curves of compact type (curves
with tree dual graphs, or equivalently, with
compact Jacobians).

Let $A^*(\overline{M}_{g})$ denote the Chow ring
with $\Q$-coefficients.
Intersection theory on $\overline{M}_g$ may be naturally viewed 
in four stages corresponding to the above filtration (\ref{fff}).
There is an associated sequence of successive
quotients:
\begin{equation}
\label{qqq}
A^*(\overline{M}_g) \rarr A^*(M^c_g) \rarr A^*(M_g) 
\rarr A^*([X_g])\eqq \Q.
\end{equation}
We develop here a uniform approach to the study of these
quotient rings.

\subsection{Tautological rings}
\label{tautor}
The study of the structure of the entire Chow ring of the
moduli space of curves appears quite
difficult at present. While presentations are known in a few
genera ([Mu], [F1], [F2], [I]), no general results have yet been
conjectured. As the principal motive is to understand 
cycle classes obtained from algebro-geometric constructions,
it is natural to restrict inquiry to the {\em tautological} ring,
$R^*(\overline{M}_g) \subset A^*(\overline{M}_g)$.

It is most convenient to define the full system of
tautological rings of all the moduli spaces of pointed
curves simultaneously:
\begin{equation}
\label{tautsys}
\{ R^*(\overline{M}_{g,n}) \subset A^*(\overline{M}_{g,n}) \}.
\end{equation}
The first step is to define the cotangent line classes
$\psi_i$. The class $$\psi_i \in A^1(\overline{M}_{g,n})$$
is the first Chern class of the line bundle with fiber
$T^*_{p_i}(C)$ over the moduli point $[C,p_1, \ldots,p_n] \in
\overline{M}_{g,n}$.
The tautological system (\ref{tautsys}) is defined
to be the set of smallest $\Q$-subalgebras  satisfying the
following three properties:
\begin{enumerate}
\item[(i)] $R^*(\overline{M}_{g,n})$ contains
the cotangent line classes $\psi_1, \ldots, \psi_n$.
\item[(ii)] The system is closed under push-forward via
all maps forgetting markings:
$$\pi_*: R^*(\overline{M}_{g,n}) \rarr R^*(\overline{M}_{g,n-1}).$$
\item[(iii)] The system is closed under push-forward via
all gluing maps:
$$\pi_*: R^*(\overline{M}_{g_1,n_1\scup\{*\}}) 
\otimes_{\Q}
R^*(\overline{M}_{g_2,n_2\scup\{\bullet\}}) \rarr
R^*(\overline{M}_{g_1+g_2, n_1+n_2}),$$
$$\pi_*: R^*(\overline{M}_{g_1, n_1\scup\{*,\bullet\}}) \rarr
R^*(\overline{M}_{g_1+1, n_1}).$$
\end{enumerate}
Natural algebraic constructions typically yield Chow classes
lying in the tautological ring.

We point out four additional properties of the
tautological system which are consequences of the above
definition:
\begin{enumerate}
\item[(iv)] The system is closed under pull-back via 
the forgetting and gluing maps.
\item[(v)] $R^*(\overline{M}_{g,n})$ is an ${\mathbb{S}}_n$-module
via the permutation action on the markings.
\item[(vi)] The $\kappa$ classes lie in the tautological rings.
\item[(vii)] The $\lambda$ classes lie in the tautological rings.
\end{enumerate}
Property (iv) follows from the well-known boundary
geometry of the moduli space of curves. As Properties
(i-iii) are symmetric under the marking permutation action,
Property (v) is obtained.
Property
(vi) is true by definition as 
$$\pi_*(\psi_{n+1}^{l+1})=\kappa_l \in R^*(\overline{M}_{g,n}),$$
where $\pi$ is the map forgetting the marking $n+1$ (see [AC]).
Recall the $\lambda$ classes are the Chern classes of the
Hodge bundle $\hodge$ on the moduli space of curves. 
Property (vii) is
a consequence of Mumford's 
Grothendieck-Riemann-Roch computation [Mu].

The tautological rings for the other elements of the
filtration (\ref{fff}) are defined by the images of $R^*(\overline{M}_g)$  
in the quotient sequence (\ref{qqq}):
\begin{equation}
\label{rrr}
R^*(\overline{M}_g) \rarr R^*(M^c_g) \rarr R^*(M_g) 
\rarr R^*([X_g])\eqq \Q.
\end{equation}

\subsection{Evaluations}
\label{eval}
The quotient rings (\ref{rrr}) exhibit several parallel
structures which serve to guide their study.
Each admits a canonical {\em non-trivial} linear
evaluation $\epsilon$ to $\Q$ obtained by integration. 
For $\overline{M}_g$, $\epsilon$ is defined
by: $$\xi\in R^*(\overline{M}_g), \ \ 
\epsilon(\xi)= \int_{\overline{M}_g} \xi.$$
The other three evaluations involve the $\lambda$ classes.

Recall the fiber of $\hodge$ over a moduli point
$[C] \in \overline{M}_g$ is the rank $g$ vector space
$H^0(C, \omega_{C})$.
Let $\Delta_0=\overline{M}_g \setminus M^c_g$. A basic
vanishing holds:
\begin{equation}
\label{kqp}
\lambda_g|_{\Delta_0}=0.
\end{equation}
To prove (\ref{kqp}), consider the standard ramified
double cover $\overline{M}_{g-1,2} \rarr \Delta_0$:
$$ [\tilde{C}, p_1,p_2] \mapsto [C]$$
obtained by identifying  the markings $p_1, p_2$ of
$\tilde{C}$ to form a nodal curve $C$.
The pull-back of $\hodge$ to $\overline{M}_{g-1,2}$
admits a surjection to the trivial bundle $\com$ over 
$\overline{M}_{g-1,2}$
obtained from the residue of $\sigma  
\in H^0(C, \omega_C)$ at the distinguished node of $C$.
Hence, the pull-back of $\lambda_g$ vanishes on $\overline{M}_{g-1,2}$.
As we consider Chow groups with $\Q$-coefficients, 
the vanishing (\ref{kqp}) follows.

For
$M^c_g$, evaluation
is defined by:
$$\xi \in R^*(M^c_g), \ \ \epsilon(\xi) =
\int_{\overline{M}_g} \xi \cdot \lambda_g,$$
well-defined by the vanishing property of $\lambda_g$.
Similarly, the vanishing of the restriction of
$\lambda_g \lambda_{g-1}$ to $\overline{M}_g \setminus
M_g$ is proven in [F3]. Define evaluation for $M_g$ by:
$$\xi \in R^*(M_g), \ \ \epsilon(\xi) =
\int_{\overline{M}_g} \xi \cdot
\lambda_g \lambda_{g-1}.$$
Finally, define evaluation for $[X_g]$ by:
$$\xi \in R^*([X_g]), \ \ \epsilon(\xi) =
\int_{\overline{M}_g} \xi \cdot \lambda_g \lambda_{g-1} \lambda_{g-2}.$$
These four evaluations do {\em not} commute with the
quotient structure.

The
non-triviality of the $\epsilon$ evaluations is proven
by explicit integral computations. 
The integral computation
\begin{equation}
\label{wwwww}
\int_{\overline{M}_{g}} \kappa_{3g-3} = \frac{1}{24^g g!}
\end{equation}
explicitly
shows $\epsilon$ is non-trivial on $R^*(\overline{M}_g)$.
Equation (\ref{wwwww}) follows from 
Witten's conjectures/Kontsevich's theorem 
or alternatively via an algebraic computation
in [FP1]. 
The integral
\begin{equation}
\label{jjj}
\int_{\overline{M}_{g}} \kappa_{2g-3} \lambda_g
 = \frac{2^{2g-1}-1}{2^{2g-1}} \frac{|B_{2g}|}{(2g)!}
\end{equation}
shows non-triviality on $R^*(M^c_g)$ [FP1]. 
The integral
\begin{equation} 
\label{sdffds}
\int_{\overline{M}_{g}} \kappa_{g-2} \lambda_g \lambda_{g-1}
= \frac{1}{2^{2g-1}(2g-1)!!} \frac{|B_{2g}|}{2g} 
\end{equation}
shows non-triviality on  $R^*(M_g)$. Equation (\ref{sdffds})
is proven in Section 1.
Finally, the computation
\begin{equation}
\label{lglglg}
\int_{\overline{M}_g} \lambda_g \lambda_{g-1} \lambda_{g-2} =
\frac{1}{2(2g-2)!}
\frac{|B_{2g-2}|}{2g-2} \frac{|B_{2g}|}{2g}
\end{equation}
establishes the last non-triviality [FP1]. 
We note the Bernoulli number convention used in these
formulas is:
\begin{equation*}
\dfrac{t}{e^t-1} = \sum\limits_{m=0}^{\infty}B_m\dfrac{t^m}{m!}.
\end{equation*}
It is known $B_{2g}$ never vanishes.

The $\epsilon$ evaluation maps are well-defined
on the
quotient sequence (\ref{qqq}) of full Chow rings.
To see difference in perspective, the non-triviality
of $\epsilon$ for $A^*(\overline{M}_g)$ is established by
considering any point class,
while the non-triviality for $R^*(\overline{M}_g)$
requires a tautological
point class --- such as a maximally degenerate stratum, or
alternatively (\ref{wwwww}).

\subsection{Gorenstein algebras}
\label{gory}
Computations of $R^*(M_g)$ for genera $g\leq 15$
have led to 
a conjecture for the ring structure
for all
genera [F3]:

\begin{conj}
$R^*(M_g)$ is a Gorenstein algebra
with socle in codimension $g-2$. 
\end{conj}

\noindent The evaluation $\epsilon$ is then
a canonically normalized function on the socle.
It is natural to hope analogous Gorenstein properties
hold for $R^*(\overline{M}_g)$ and $R^*(M^c_g)$, but
the data in these cases is very limited. The following
conjectures are therefore really speculations.

\vspace{+10pt}
\noindent {\bf Speculation 2.} {\em 
$R^*(M^c_g)$ is a Gorenstein algebra with socle
in codimension $2g-3$.}

\vspace{+10pt}
\noindent {\bf Speculation 3.} {\em
$R^*(\overline{M}_g)$ is a Gorenstein algebra
with socle in codimension $3g-3$.}

\vspace{+10pt}

Conjecture 1 was verified for $g\leq 15$ via
relations found by classical degeneracy loci techniques [F3]
and the non-vanishing result (\ref{sdffds}) --- see Section 1.
In fact, a complete presentation of $R^*(M_g)$
has been conjectured in [F3] from these low genus studies.
Such calculations become much more difficult in
$R^*(M_g^c)$ and $R^*(\overline{M}_g)$ because of
the inclusion of nodal curves. $R^*(M_g^c)$ and $R^*(\overline{M}_g)$
are known to be Gorenstein algebras for $g\leq 3$.
It would be very interesting to find further evidence
for or against Speculations 2 and 3. 

A stronger version of Conjecture 1 was made in [HL].
Also, Speculation 3 was raised as a question in [HL].

An extension of the perspective presented here to pointed curves
and fiber products of the universal curve will be discussed
in [FP3].

As the moduli space of stable curves  $\overline{M}_{g,n}$ may be viewed
as a special case of the moduli
space of stable maps $\overline{M}_{g,n}(X,\beta)$, it is natural
to investigate tautological rings in the more general setting
of stable maps.
The first obstacle is finding the appropriate definitions
in the context of the virtual class. However, in the case
of genus 0 maps to homogeneous varieties, it is straightforward
to define the tautological ring since the moduli space is
a nonsingular Deligne-Mumford stack. 
In [P1], the tautological ring
$R^*(M_{0,0}(\proj^r,d))$ is proven to be a Gorenstein algebra.

\subsection{Socle rank and  higher vanishing predictions}
The Gorenstein Conjectures/ Speculations of Section \ref{gory}
imply the ranks  of the tautological rings
are 1 in the expected socle codimension.
Moreover, vanishing above the socle codimension is implied in each case.
The socle and vanishing results
 $$R^{g-2}(M_g)\eqq \Q, \ \ R^{>g-2}(M_g)=0 $$
are a direct
consequence of 
Looijenga's Theorem [L] and the non-vanishing (\ref{sdffds}) proven
in Section 1.
Looijenga's Theorem states the tautological ring of the $n$-fold
fiber product $C_g^n$ of $C_g=M_{g,1}$ over $M_g$ 
is {\em at most} rank 1 in codimension $g-2+n$ and 
vanishes in all codimensions
greater than
$g-2+n$.

It is natural to ask whether the tautological rings
satisfy the usual right exact sequences via restriction:
\begin{equation}
\label{exexexex}
R^*(\partial \overline{M}_g) \rarr R^*(\overline{M}_g)
\rarr R^*(M_g) \rarr 0.
\end{equation}
Here, $R^*(\partial 
\overline{M}_g)\subset A^*(\partial \overline{M}_g)$ is 
generated by 
tautological classes pushed forward to the boundary
$\partial \overline{M}_g$ of the moduli space of curves.
Pointed generalizations of the restriction
 sequences (\ref{exexexex}) together with 
Looijenga's Theorem and the non-vanishings (\ref{wwwww}-\ref{jjj})
imply the socle and vanishing results for 
$R^*(M_g^c)$ and $R^*(\overline{M}_g)$. However, at present, the
right exactness of sequence
(\ref{exexexex}) is not proven. 

We note the socle dimension
proof for $R^*(\overline{M}_g)$ in Section 5.1 of [HL]
is incomplete as it stands since (\ref{exexexex}) is assumed there
(the error is repeated in [FL]).

\subsection{Virasoro constraints}
\label{wwww}
The tautological rings  (\ref{rrr}) each have an associated Virasoro
conjecture. For $\overline{M}_g$, the original Virasoro
constraints (conjectured by Witten and proven by Kontsevich [K1])
compute all the
integrals
\begin{equation}
\label{zzx}
\int_{\overline{M}_{g,n}} \psi_1^{\alpha_1} \cdots \psi_n^{\alpha_n}.
\end{equation}
These integrals determine the $\epsilon$ 
evaluations in the ring $R^*(\overline{M}_g)$. The methods
for calculating $\epsilon$ evaluations 
from the integrals (\ref{zzx}) are effective but quite complicated
(see [F3], [HL], [W]).

Eguchi, Hori, and Xiong (and S. Katz)
have conjectured Virasoro constraints
in Gromov-Witten theory for general target varieties $V$
which specialize to Witten's conjectures in case $V$ is a point
[EHX]. In [GeP],
these general constraints are applied
to {\em collapsed} maps to target  curves, surfaces, and threefolds
in order to study integrals of the Chern classes of the Hodge bundle.
The Virasoro constraints for curves then 
imply: 
\begin{equation} \label{lamg}
\int_{\overline{M}_{g,n}} 
\psi_1^{\alpha_1} \cdots \psi_n^{\alpha_n} \lambda_g =
\binom{2g+n-3}{\alpha_1,\ldots, \alpha_n}\int_{{\overline{M}}_{g,1}}
\psi_1^{2g-2} \lambda_g,
\end{equation}
where $\alpha_i\geq 0$.
Equation (\ref{lamg}) determines (up to scalars)
the $\epsilon$ evaluations  in the ring $R^*(M^c_g)$.
This Virasoro conjecture for $M_g^c$ has been proven in
[FP2].

The Virasoro constraints for surfaces imply a formula
previously conjectured in [F3] determining
evaluations in $R^*(M_g)$:
\begin{equation} \label{lamgg}
\int_{\overline{M}_{g,n}} 
\psi_1^{\alpha_1} \cdots \psi_n^{\alpha_n} \lambda_g \lambda_{g-1}
= \frac{(2g+n-3)! (2g-1)!!}{(2g-1)!\prod_{i=1}^n (2\alpha_i-1)!!}
\int_{\overline{M}_{g,1}} \psi_1^{g-1} \lambda_g \lambda_{g-1},
\end{equation}
where $\alpha_i>0$ (see [GeP]).
Formula (\ref{lamgg}) is currently still conjectural.

Finally, the Virasoro constraints for threefolds yield relations among
the integrals
\begin{equation} \label{lamggg}
\int_{\overline{M}_{g,n}} 
\psi_1^{\alpha_1} \cdots \psi_n^{\alpha_n} \lambda_g \lambda_{g-1}
\lambda_{g-2}.
\end{equation}
In fact, all integrals (\ref{lamggg}) are determined 
in terms of $\int_{\overline{M}_{g}} \lambda_g \lambda_{g-1} \lambda_{g-2}$
by the string and dilaton equations (which leads to a proof
of the Virasoro constraints in this case [Ge]). 

We note the ring structure of a
finite dimension Gorenstein algebra is {\em determined}
by the socle evaluation of polynomials in the generators.
Hence, 
if the Gorenstein properties of Section \ref{gory} hold
for any of the tautological rings, the Virasoro constraints
then determine the ring structure.
This concludes our general discussion of the tautological
rings of the moduli space of curves.

\subsection{Results}
\label{sss}
A basic generating series for $1$-pointed 
Hodge integrals was computed in [FP1]:
\begin{equation}
\label{sdsd}
1+ \sum_{g\geq 1} \sum_{i=0}^{g} t^{2g} k^i
\int_{\overline{M}_{g,1}} \psi_1^{2g-2+i} \lambda_{g-i} = 
\Big( \frac{t/2}{\sin(t/2)} \Big)^{k+1}.
\end{equation}
Equation (\ref{sdsd})
may be interpreted to determine $\epsilon$ evaluations of the monomials 
$$\kappa_{3g-3-i} \lambda_i \in
R^{3g-3}(\overline{M}_{g}).$$
The main result of this paper is a determination of related evaluations
in $R^{g-2}(M_g)$.

First, the basic series for the non-triviality of $\epsilon$
on $R^*(M_g)$ is calculated.

\vspace{+10pt}
\noindent{\bf Theorem 1.} {\em 
For genus $g\geq 2$},
\begin{equation}
\label{qwqwz}
\int_{\overline{M}_{g}} \kappa_{g-2} \lambda_g \lambda_{g-1}
= \frac{1}{2^{2g-1}(2g-1)!!} \frac{|B_{2g}|}{2g}. 
\end{equation}
\vspace{+10pt}

\noindent 
Two proofs of Theorem 1 are given in the paper.
The first 
uses Mumford's  Grothendieck-Riemann-Roch 
formulas for the Chern character of $\hodge$ and
the Witten/Kontsevich theorem
in KdV form. 
The derivation appears in Section 1, following a discussion of the
context of this calculation.
The second proof appears in Section 5 as a combinatorial
consequence of 
Theorem 3 below. The required combinatorics is explained
in the Appendix by D. Zagier.

Next, integrals  encoding the values of all the monomials
$$ \kappa_{g-2-i} \lambda_i \in R^{g-2}(M_g)$$
are studied.
For positive integers $g$ and $k$, let
$$I(g,k)= \int_{{\overline{M}_{g,1}}}
\frac{ 1-\lambda_1 +\lambda_2-\ldots + 
(-1)^g \lambda_g}{ \prod_{i=1}^{k} (1-i\psi_1)} 
\ \lambda_g \lambda_{g-1}.$$
The integrals $I(g,k)$ arise geometrically in 
the following manner. Let $$\pi: M_{g,1} \rarr M_g$$
be the universal curve. Let $J_{k}$ denote
the rank $k$ vector bundle with fiber 
$$H^0(C, \omega_C/\omega_{C}(-kp))$$
at the moduli point $[C,p]$. $J_{k}$ is a bundle of
$\pi$-vertical $(k-1)$-jets of $\omega_\pi$. 
There is a canonical (dualized) evaluation map
\begin{equation}
\label{twq}
J_{k}^* \rarr \hodge^*
\end{equation}
on $M_{g,1}$. For $g\geq 2$, 
$$I(g,k)= \epsilon( \  \pi_* c_{g-1}\Big(\frac{\hodge^*}{J^*_{k}} \Big)\ 
),
$$
where the $\epsilon$ evaluation is taken in $R^*(M_g)$.

For $k=1$, $J_1= \omega_\pi$ and
 the map (\ref{twq}) is a bundle injection.
$I(g,1)$ is then the evaluation of the $\pi$-push forward
of the Euler class of the quotient:
$$I(g,1)= \epsilon( \  \pi_* c_{g-1}\Big(\frac{\hodge^*}{\omega_\pi^*} 
\Big)\ 
).
$$
The integrals $I(g,2)$ are easily related to the
(stack) classes of the hyperelliptic loci $[H_g] \in R^{g-2}(M_g)$
by the equation (see [Mu]):
\begin{equation}
\label{mwm}
I(g,2)= (2g+2)\cdot \epsilon([H_g]).
\end{equation}
For $k >2$, $I(g,k)$ does not admit such simple
interpretations. However, generating series
of these integrals appear to be the best behaved
analogues of (\ref{sdsd}) in $R^*(M_g)$. The search
for such an analogue was motivated by the 
parallel structure view of these tautological rings.

For each positive integer $k$, define
$$G_k(t)=  \sum_{g \geq 1} t^{2g+k-1}  I(g,k).$$
These generating series are uniquely determined by:

\vspace{+10pt}
\noindent{\bf Theorem 2.} {\em 
For all integers $k\geq 1$,
the series $G_k(t)$ satisfies}
\begin{equation}
\label{qwqw}
 \frac{d^{k-1} G_k}{dt^{k-1}} 
= \sum_{j=1}^k { (-1)^{k-j}}
\frac{j^{k-1}}
{k} \binom {k}{j} 
 {\text {log}} \Big( \frac{jt/2}{{\text{sin}}(jt/2)} \Big).
\end{equation}
\vspace{+10pt}

In case $k=1$, we obtain the following Corollary
first encountered in the study of degenerate 3-fold
contributions in Gromov-Witten theory [P2].

\vspace{+10pt}
\noindent{\bf Corollary 1.} 
$$\sum_{g \geq 1} t^{2g} \int_{\overline{M}_g}
 (\sum_{i=0}^{g-2} (-1)^i \kappa_{g-2-i} \lambda_i)\ \lambda_g 
\lambda_{g-1}
= {\text {log}} \Big( \frac{t/2}{{\text{sin}}(t/2)} \Big).$$
\vspace{+10pt}

In case $k=2$, we find
$$(G_2)'= {\text{log}} \Big( \frac{2\text{sin}(t/2)}{\text{sin}(t)} 
\Big)
= - {\text{log}}(\text{cos}(t/2)).$$
The generating series for the evaluations of the
hyperelliptic loci in $R^*(M_g)$ (with
an appropriate genus 1 term) is:
$$H(t)= \frac{t^2}{96} + \sum_{g\geq2} t^{2g} \epsilon([H_g]).$$
By Mumford's calculation (\ref{mwm}),
\begin{equation*}
 (t^2 H)'= G_2
\end{equation*}
Theorem 2 then yields the following result.

\vspace{+10pt}
\noindent{\bf Corollary 2.} 
{\em The hyperelliptic evaluations are determined by:}
\begin{equation}
\label{zzz}
(t^2 H)''= -{\text{log}}({\text{cos}}(t/2)).
\end{equation}
\vspace{+10pt}

\noindent Equation (\ref{zzz}) was  conjectured previously in an
equivalent Bernoulli number form in [F3]: for $g\geq 2$,
$$\epsilon([H_g])= 
\frac{ (2^{2g}-1)|B_{2g}|}{(2g+2)! \ 2g}.
$$

Theorem 2 is derived here from relations obtained
by virtual localization in Gromov-Witten theory (see [GrP], [FP1],
[FP2]).
In addition to the cohomology classes on the moduli
space of stable maps $\overline{M}_{g,n}(\proj^1,d)$
considered in [FP2], new classes obtained
from the ramification map of [FanP] play
an essential role. The Hodge integral series (\ref{sdsd}) and
Virasoro constraints (\ref{lamg}) for
$M_g^c$ are also used.
This derivation appears in
Sections 2 and 3 of the paper.

In case
$k=2$, the integrals $I(g,2)$
may be computed by reduction to the moduli space
of hyperelliptic curves.
This classical derivation 
provides a  contrast to the more formal
Gromov-Witten arguments. 
Section 4 of the paper contains these hyperelliptic
computations.

In Section 5, the standard 1-point Hodge integral series for
$R^*(M_g)$ is studied. The following consequence of Theorem 
2 is found.

\vspace{+10pt}
\noindent{\bf Theorem 3.}
{\em
For positive integers $g$, $k$,}

$$
\sum_{i=0}^{g-1}(-1)^ik^{g-1-i}
\int_{\overline{M}_{g,1}} \psi_1^{g-1-i}\lmb{i}\lmb{g}\lmb{g-1} 
=
\frac{|B_{2g}|}{2g}\,
\sum_{l=1}^k
\frac{(k-1)!}{(k-l)!}\,\frac{l!}{k^l} \frac{{\mathfrak {S}}^{(l)}_{2g-1+l}}
{(2g-1+l)!}\,. 
$$
\vspace{+10pt}

\noindent
Here, $\mathfrak{S}^{(l)}_{n+l}$ is the Stirling number of the
second kind: $\mathfrak{S}^{(l)}_{n+l}$ equals the number of partitions
of a set of $n+l$ elements into $l$ non-empty subsets.

Theorem 3 and the Appendix together provide
proofs of all previously conjectured formulas 
for 1-point integrals 
in the tautological ring. In particular, closed forms for
the evaluations in $R^*(M_g)$
of
\begin{equation}
\label{werrew}
\kappa_{g-2},\  \kappa_{g-3} \lambda_1, \ \kappa_{1}\lambda_{g-3},
\ \lambda_{g-2}
\end{equation}
are found --- providing an alternate derivation of Theorem 1
and settling conjectures of [F3],[F4].
A list of these formulas
is provided in Section 5.2. 
In fact, the combinatorial
results of the Appendix lead to 
proofs of natural extensions of the formulas for
(\ref{werrew}).

\subsection{Acknowledgements}
We thank  D. Zagier for his aid in our work --- especially
for the results proven in the Appendix. Also, conversations
with R. Dijkgraaf, E. Getzler, and S. Popescu were helpful to us. 
The authors were partially supported by 
National Science
Foundation grants DMS-9801257 and DMS-9801574.
C.~F.~was partially supported by the Max-Planck-Institut
f\"ur Mathematik, Bonn.
R.~P.~was partially supported by an A.~P.~Sloan foundation
fellowship.

\section{\bf{Theorem 1}}

\subsection{Context} 
Looijenga has proven in [L] that the tautological ring $R^*(M_g)$ vanishes
in degrees greater than $g-2$ and is at most one-dimensional in 
degree $g-2$, generated by the class of the hyperelliptic locus.
Theorem 1 
shows 
$$\epsilon(\kp{g-2}) = \int_{\overline{M}_g}
\kp{g-2}\lmb{g}\lmb{g-1}$$ is nonzero where $\epsilon$
is the evaluation on $R^*(M_g)$, see Section \ref{eval}. Hence,
$\kp{g-2}$ is nonzero in $R^{g-2}(M_g)$.
In Section 1.2, we present the first proof (Fall 1995) of Theorem 1
--- relying upon
an explicit calculation using the Witten/Kontsevich theorem
in KdV form. The resulting non-vanishing of the tautological ring $R^*(M_g)$
in degree $g-2$ completed the verification
for $5\le g\le 15$ of the
conjectural description of $R^*(M_g)$ given in [F3].
A second, more geometric proof of this non-vanishing appears in Section 4
using the defining property of hyperelliptic curves.
Later proofs may be found in [FP1] and [P2], showing the non-vanishing
in $R^{g-2}(M_g)$
of $\lmb{g-2}$ and $\sum_{i=0}^{g-2}(-1)^i\kp{i}\lmb{g-2-i}$
respectively.
Theorem 1 is rederived in Section 5 from Theorem 3 (together with the
Appendix) providing an alternative
to the KdV derivation here.

\subsection{First proof of Theorem 1}
Using Mumford's expression [Mu] for the Chern character of the Hodge bundle
and the resulting identity [FP1]
$$\lmb{g}\lmb{g-1}=(-1)^{g-1}(2g-1)!\,\ch_{2g-1}(\hodge),$$
Theorem 1 is reduced to the identity
\begin{equation}
\label{nv1}
\frac1{2^{2g-1}(2g-1)!!} =
\lan\ta{2g}\ta{g-1}\ran-\lan\ta{3g-2}\ran+{\frac12}\sum_{j=0}^{2g-2}
(-1)^j\lan\ta{2g-2-j}\ta{j}\ta{g-1}\ran\hfill
\end{equation}
$$
{}+{\frac12}\sum_{h=1}^{g-1}\left((-1)^{g-h}\lan\ta{3h-g}\ta{g-1}\ran
\lan\ta{3g-3h-2}\ran
+(-1)^h\lan\ta{3h-2}\ran\lan\ta{2g-3h}\ta{g-1}\ran\right) $$
(see [FP1]). Here, the second sum equals
$$
\sum_{h=1}^{g-1}{\frac{(-1)^{g-h}}{24^{g-h}(g-h)!}}\lan\ta{3h-g}\ta{g-1}\ran
$$
since $\lan\ta{3k-2}\ran=1/(24^kk!)$ by equation (\ref{wwww}).
Hence, it suffices to prove the two identities
\begin{equation}
\label{four}
\sum_{h=1}^{g}{\frac{(-1)^{g-h}}{24^{g-h}(g-h)!}}\lan\ta{3h-g}\ta{g-1}\ran
={\frac1{24^gg!}}
\end{equation}
and
\begin{equation}
\label{five}
\sum_{j=0}^{2g-2}(-1)^j\lan\ta{2g-2-j}\ta{j}\ta{g-1}\ran=
{\frac{g!}{2^{g-2}(2g)!}}\quad.
\end{equation}
Both are consequences of the following  equation for coefficients resulting
from Witten's KdV-equation for power series ([W], (2.33), (2.19)). 
For any monomial
$$ T = \prod_{j=0}^k \tau_j^{d_j}, $$
the coefficient  equation holds:
\begin{equation}
\label{brief0}
(2n+1)\lan\tau_n\tau_0^2 T\ran = {\frac14}\lan\tau_{n-1}\tau_0^4T\ran 
\end{equation}
$$
+\sum_{0\le a_j\le d_j}
\left(\prod_{j=0}^k {\binom{d_j}{a_j}}\right)
\bigl(\lan\tau_{n-1}\tau_0T_1\ran\lan\tau_0^3T_2\ran +
 2\lan\tau_{n-1}\tau_0^2T_1\ran\lan\tau_0^2T_2\ran\bigr)    
$$
where the sum is over factorizations $T=T_1T_2$ with 
$T_1=\prod_{j=0}^k \tau_j^{a_j} $.

For $T=\tau_b$ and $n=a$, this gives
\begin{equation}
\label{brief1}
(2a+1)\lan\tau_0^2\tau_a\tau_b\ran ={\frac14}\lan\tau_{a-1}\tau_0^4\tau_b\ran
+\lan\tau_{a-1}\tau_0\tau_b\ran\lan\tau_0^3\ran
\end{equation} $$
\qquad{}+
\lan\tau_{a-1}\tau_0\ran
\lan\tau_0^3\tau_b\ran+2\lan\tau_{a-1}\tau_0^2\tau_b\ran\lan\tau_0^2\ran
+2\lan\tau_{a-1}\tau_0^2\ran\lan\tau_0^2\tau_b\ran.
$$
Consider now the two-point function 
$D(w,z)=\sum_{a,b\ge0}\lan\ta0\ta{a}\ta{b}\ran w^az^b$.
Equation (\ref{brief1}) is equivalent to the differential equation:
\begin{equation}
\left(2w{\frac{\pa}{\pa w}}+1\right)\bigl((w+z)D(w,z)\bigr)
={\frac14}(w+z)^3wD(w,z)+wD(w,z)
\end{equation}
$$
\qquad{}+D(w,0)zD(0,z)+2wD(w,0)D(0,z)\,.
$$
It is easy to verify that the unique solution of this equation satisfying
$D(w,0)=\exp(w^3/24)$ and $D(0,z)=\exp(z^3/24)$ is given by
$$
D(w,z)=\exp\left({\frac{(w^3+z^3)}{24}}\right)
\sum_{n\ge0}{\frac{n!}{(2n+1)!}}
\left[{\textstyle{\frac12}}wz(w+z)\right]^n\,.
$$
We learned this formula from Dijkgraaf [Dij].
Consequently, for all $k\ge1$
\begin{equation}
\label{dfgf}
\sum_{h=0}^g{\frac{(-1)^{g-h}}{24^{g-h}(g-h)!}}\lan\ta{0}\ta{3h-g+k}\ta{g-k}\ran
=0\,,
\end{equation}
since this is the coefficient of $w^{2g+k}z^{g-k}$ in
$$
\lan\ta{0}\tau(w)\tau(z)\ran\cdot\lan\ta{0}\tau(-w)\ta{0}\ran
=\exp\left({\frac{z^3}{24}}\right)\sum_{n\ge0}{\frac{n!}{(2n+1)!}}
\left[{\textstyle{\frac12}}wz(w+z)\right]^n\,,
$$
in which all terms of total degree $3g$ have degree at least $g$
in $z$. Therefore, by applications of the string equation to (\ref{dfgf}),
we find:
\begin{eqnarray*}
&&\sum_{h=0}^g{\frac{(-1)^{g-h}}{24^{g-h}(g-h)!}}\lan\ta{3h-g}\ta{g-1}\ran
=-\sum_{h=0}^g{\frac{(-1)^{g-h}}{24^{g-h}(g-h)!}}\lan\ta{3h-g+1}\ta{g-2}\ran\cr
&&=
+\sum_{h=0}^g{\frac{(-1)^{g-h}}{24^{g-h}(g-h)!}}\lan\ta{3h-g+2}\ta{g-3}\ran
=\dots          \cr
&&=(-1)^{g-1}
\sum_{h=0}^g{\frac{(-1)^{g-h}}{24^{g-h}(g-h)!}}\lan\ta{3h-1}\ta{0}\ran
=\sum_{h=1}^g{\frac1{24^hh!}}\,{\frac{(-1)^{h+1}}{24^{g-h}(g-h)!}}  \cr
&&={\frac1{24^gg!}}\left(\sum_{h=1}^g(-1)^{h+1}{\binom{g}h}\right)=
{\frac1{24^gg!}}\quad,\cr
\end{eqnarray*}
which proves (\ref{four}) for $g\ge1$.

To prove (\ref{five}), we use (\ref{brief0}) for
$T=\ta{b}\ta{c}$ and $n=a$. This is equivalent 
to a differential equation for the general three-point function
$E(x,y,z)=\sum_{a,b,c\ge0}\lan\ta{a}\ta{b}\ta{c}\ran x^ay^bz^c$
that specializes to the following differential equation
for the special three-point function
$F(w,z)=E(w,z,-z)$:
$$\displaylines{\quad
4w^2F(w,z)+2w^3{\frac{\pa F}{\pa w}}(w,z)-{\frac14}w^5F(w,z)\hfill\cr
\hfill{}=w(2w+z)D(w,z)D(0,-z)+w(2w-z)D(w,-z)D(0,z)\quad.\quad\cr}
$$
It is clear that it has a unique solution. One verifies easily
that the solution is
$$
F(w,z)=\exp\left({\frac{w^3}{24}}\right)\sum_{a,b\ge0}
(w^3)^a(wz^2)^b{\frac{(a+b)!}{2^{a+b-1}(2a+2b+2)!}}
{\binom{a+b+1}{2a+1}}\quad.
$$
The coefficient of $w^gz^{2g}$ equals
$${\frac{(g+1)!}{2^{g-1}(2g+2)!}}\quad,
$$
which gives (\ref{five}).
This finishes the (first) proof of Theorem 1.

\section{\bf{Localization relations}}

\subsection{Results}
In this Section, the localization method 
will be used to find relations among Hodge integrals [FP1], [FP2].
Define the Hodge integral $Q^e_g$ for $g,e\geq 1$ by:
\begin{equation}
\label{mqz}
Q^e_g= \int_{\overline{M}_{g,1}} 
\frac{ 1-\lambda_1 +\lambda_2-\ldots + 
(-1)^g \lambda_g}{ 1-e \psi_1 } 
\ \lambda_g \lambda_{g-1}.
\end{equation}
The first step in the proof of Theorem 2
 is the computation of $Q^e_g$.

To state the relations determining $Q^e_g$, we will need the
following combinatorial coefficients.
For any formal series $t(x)=\sum {t_i} x^i$ define
$${\mathcal{C}}(x^i,t(x))= t_i.$$
Let $\tau(x)$ be the series inverse of $xe^{-x}$:
$$\tau(x)=\sum_{r\geq 1} \frac{r^{r-1}}{r!} x^r.$$
For $d\geq e$, define $f_{gde}$ by:
\begin{equation}
\label{jack}
f_{gde}= 
\frac{e^{e+1}}{e!} \sum_{l=0}^{2g} \frac{(2g+d-l-1)!}{(2g-l)!}
\frac{(-d)^l}{l!}
 \ {{\mathcal C}( x^{d-e}, \tau^l(x) )}.
\end{equation}

\begin{pr} \label{rell}
For  $d\geq 1$,
$$\sum_{e=1}^d \sum_{g=1}^\infty Q_g^e  f_{gde} t^{2g}
= d^{d-1} 
{\text {\em log}}
\Big( \frac{dt/2}{{\text{\em sin}}(dt/2)} \Big).$$
\end{pr}

\noindent
The proof of Proposition \ref{rell} depends upon almost all of the
main results of [GrP], [FP1], [FP2], and [FanP].
Theorem 2 will be derived as a consequence of Proposition \ref{rell}
in
Section 3.

\subsection{The torus action}
Let $\proj^1=\proj(V)$ where $V=\com \oplus \com$.
Let $\com^*$ act diagonally on $V$:
\begin{equation}
\label{repp}
\xi\cdot (v_1,v_2) = ( v_1, 
\xi \cdot v_2).
\end{equation}
Let $p_1, p_2$ be the fixed points $[1,0], [0,1]$ of the corresponding
action on $\proj(V)$.
An equivariant lifting  of $\com^*$ to a line bundle $L$
over 
$\proj(V)$ is uniquely determined by the weights $[l_1,l_2]$
of the fiber
representations at the fixed points 
$$L_1= L|_{p_1}, \ \ \ 
L_2= L|_{p_2}.$$
The canonical lifting of $\com^*$ to the
tangent bundle $T_\proj$ has weights $[1,-1]$.
We will utilize the equivariant liftings of
$\com^*$ to $\oh_{\proj(V)}(1)$ and $\oh_{\proj(V)}(-1)$ with weights
$[0,-1]$, $[0,1]$ respectively.

Let $\overline{M}_{g,n}(\proj(V), d)$ be the moduli
stack of stable genus $g$, degree $d$ maps to $\proj^1$ (see [K2], [FuP]).
There are canonical maps
$$\pi: U \rarr \overline{M}_{g,n}(\proj(V),d), 
\ \ \ \mu: U \rarr \proj(V)$$
where $U$ is the universal curve over the moduli stack.
The representation (\ref{repp}) canonically
induces $\com^*$-actions on $U$ and 
$\overline{M}_{g,n}(\proj(V),d)$ compatible
with the maps $\pi$ and $\mu$ (see [GrP]).

\subsection{The branch morphism}
In [FanP], a canonical branch divisor morphism
$\gamma$ is constructed using derived category techniques:
\begin{equation}
\label{exxx}
\gamma: \overline{M}_{g,n}(\proj(V),d)\rarr  \sym^{r}(\proj(V))=
\proj(\sym^r(V^*)),
\end{equation} 
where $r=2d+2g-2$.
We review the 
point theoretic
description of $\gamma$.
Let 
$$[f: C \rarr \proj(V)]$$ be a moduli point where
$C$ is a possibly singular curve. Let $N\subset C$ be the
cycle of nodes of $C$.
Let $\nu: \tilde{C} \rarr C$
be the normalization of $C$.
Let $A_1, \ldots, A_a$ be the components of $\tilde{C}$
which dominate $D$, and let $$\{a_i: A_i \rarr D\}$$
denote the natural maps. 
As $a_i$ is a surjective map between nonsingular
curves, the classical branch divisor $\br(a_i)$ is well-defined. 
Let $B_1, \ldots, B_b$ be
the components of $\tilde{C}$ contracted over $D$, and
let $f(B_j)=p_j\in D$.
Then, the following formula holds:
\begin{equation}
\label{ptwise}
\gamma([f])=\br(f)= \sum_{i} \br(a_i) +
\sum_j (2g(B_j)-2)[p_j] + 2f_{*}(N).
\end{equation}
We note $\gamma$ commutes with the forgetful maps
$$\overline{M}_{g,n}(\proj(V),d) \rarr \overline{M}_g(\proj(V),d),$$
and $\gamma$ is equivariant with respect to the
canonical action of $\com^*$ defined by the 
representation (\ref{repp}).

\subsection{Equivariant cycle classes}
\label{ecc}
We now describe the  equivariant Chow classes which
arise in the proof of Proposition 1.

First consider the $\com^*$-action on $\proj(\sym^r(V^*))$.
There are exactly $r+1$ distinct
$\com^*$-fixed points.
For $0\leq a \leq r$,
let $q_{a}$ denote the fixed point $v_1^{*(r-a)} v_2^{*a}$. 
The canonical $\com^*$-linearization on $S= \oh(1)$
has weight $$w_a=a$$ at $q_{a}$.
Let $S_i$ denote the unique 
$\com^*$-linearization of $S$
satisfying $w_{i}=0.$
We note the weight at $q_a$ of $S_i$ is $a-i$. 
The first equivariant Chow classes considered are
$$s_i= \gamma^*(c_1(S_i)),$$
for all $0\leq i \leq r$.

Second,
there is a natural rank $d+g-1$ bundle on 
$\overline{M}_{g,n}(\proj(V),d)$: 
\begin{equation}
\label{wqwq}
\mathbb{R}=R^1\pi_* (\mu^* \oh_{\proj(V)}(-1)).
\end{equation}
The linearization $[0,1]$ on $\oh_{\proj(V)}(-1)$
defines an equivariant $\com^*$-action on $\mathbb{R}$.
We will require the equivariant top Chern class $c_{top}(\mathbb{R})$.

Third,
there is a canonical lifting of the 
$\com^*$-action on $\overline{M}_{g,n}(\proj(V),d)$ to the
Hodge bundle $\hodge$ over $\overline{M}_{g,n}(\proj(V),d)$.
Hence, the Chern classes $\lambda_i$ yield equivariant
cycle classes.

Finally, let 
$$\text{ev}_i: \overline{M}_{g,n}(\proj(V),d) \rarr \proj(V)$$
denote the $i^{th}$ evaluation morphism, and let
$$\rho_i =c_1( \text{ev}_i^* \oh_{\proj(V)}(1)),$$
where we fix the $\com^*$-linearization $[0,-1]$ on
$\oh_{\proj(V)}(1)$.

\subsection{Vanishing integrals}
We will obtain relations among $Q^e_g$ from
a sequence of vanishing integrals.
Let $g,d \geq 1$.
Let $P(g,d)$ denote the integral:
$$P(g,d)= \int_{\overline{M}_{g,1}(\proj^1,d)}
\lambda_{g-1}\  c_{top}(\mathbb{R})\ \rho_1^2 \ \prod_{i=0}^{d-2} s_i =0.$$
As the virtual dimension of $\overline{M}_{g,1}(\proj^1,d)$
equals $2d+2g-1$ and the total dimension of the integrand
is $$(g-1)+(d+g-1)+2+(d-1)=2d+2g-1,$$
the integral $P(g,d)$ is well-defined.
Since $\rho_1^2=0$, $P(g,d)=0$.

\subsection{Localization terms}
\label{lres}
As all the integrand terms in $P(g,d)$
have been defined with $\com^*$-equivariant lifts,
the virtual localization formula of [GrP] yields
a computation of these integrals in terms
of Hodge integrals over moduli spaces of stable curves.

The integrals
$P(g,d)$ are
expressed via localization as a  sum over connected decorated 
graphs $\Gamma$ (see [K2], [GrP]) indexing the $\com^*$-fixed loci of 
$\overline{M}_{g,n}(\proj(V),d)$.
The vertices of these graphs lie over the
fixed points $p_1, p_2 \in \proj(V)$ and are
labelled with genera (which sum over the graph to $g-h^1(\Gamma)$).
The edges of the graphs lie over $\proj^1$ and
are labelled with degrees (which sum over the 
graph to $d$). Finally, the graphs carry
a single marking on one of the vertices. The edge valence of a vertex
is the number of incident edges (markings excluded).

The equivariant integrand of $P(g,d)$
has been chosen to force
vanishing contributions for most graphs (see [FP1], [FP2]).
By the linearization choice 
on the bundle $\mathbb{R}$, we find:
if a graph $\Gamma$
contains a vertex lying over $p_1$ of 
edge valence greater than 1, then the contribution of
$\Gamma$ to $P(g,d)$ vanishes.
This basic vanishing was first used in $g=0$ by Manin in [Ma].
Additional applications have been pursued in [GrP], [FP1], [FP2].

By the above vanishing, only {\em comb} graphs $\Gamma$ 
contribute to $P(g,d)$. Comb graphs contain
$k\leq d$ vertices lying over $p_1$ each
connected by a distinct edge to a unique vertex lying
over $p_2$. These graphs carry the usual vertex genus
and marking data. 

If the (unique) marking of $\Gamma$ lies over $p_1$, then
the contribution of $\Gamma$ to $P(g,d)$ vanishes by the
linearization choice for $\rho_1$. We may thus assume the
marking of $\Gamma$ lies over $p_2$. 

A comb graph $\Gamma$ is defined to have complexity $n\geq 0$ if
exactly $n$ vertices lying over $p_1$ have positive genus.
A vertex $v$ of positive genus $g(v)$ over $p_1$ 
yields the moduli space $\overline{M}_{g(v),1}$
occurring as a factor in the fixed point locus
corresponding to $\Gamma$. Let $v_1, \ldots, v_{k'}$
denote the positive genus vertices over $p_1$.
The fixed point locus corresponding to $\Gamma$ 
is a quotient of
\begin{equation}
\label{lqrp}
\prod_{i=1}^{k'} \overline{M}_{g(v_i),1} 
\times \overline{M}_{g',k+1}.
\end{equation}
Here, the unique vertex over $p_2$
is of genus $g'$, the comb consists of $k$ total
vertices over $p_1$, and the marking lies over $p_2$.
The restriction of the integrand term $c_{top}({\mathbb R})$
to the fixed locus yields the class  
$$\prod_{i=1}^{k'} \lambda_{g(v_i)}$$
as a factor.
The integrand term $\lambda_{g-1}$ contributes
the sum:
\begin{equation}
\label{sdhj}
\prod_{i=1}^{k'} \lambda_{g(v_i)} \ \lambda_{g'-1}\ +
\sum_{i=1}^{k'} \lambda_{g(v_i)-1} \prod_{j\neq i} \lambda_{g(v_j)}\
\lambda_{g'}.
\end{equation}
By (\ref{sdhj}) and the basic vanishing
$\lambda_h^2=0 \in A^*(\overline{M}_{h,1})$ for $h>0$,
we easily see
graphs $\Gamma$ of complexity greater than 1 contribute
0 to $P(g,d)$.
We have proven only graphs of complexity 0 or 1 may contribute
to $P(g,d)$.

Consider first a graph $\Gamma$ of complexity 0. As before,
let $k$ be the total number of vertices over
$p_1$.
The image under $\gamma$ 
of the fixed point locus corresponding to $\Gamma$ is the point
$q_{d-k}$. By the term $\prod_{i=0}^{d-2} s_i$ in the integrand,
all such graphs contribute 0 unless $k=1$.
Therefore there is a unique complexity 0 graph $\Gamma$ which
contributes to $P(g,d)$.
The contribution of this graph is:
\begin{equation}
\label{jjjj}
-(-1)^{d-g} \ d^{d-2} \ d^{2g} \int_{\overline{M}_{g,1}} \psi_1^{2g-1}
\lambda_{g-1}.
\end{equation}  
The contribution is computed via a direct application of
the virtual localization formula [GrP]. Only one Hodge integral
(occurring at the vertex lying over $p_2$) appears.
 
Next, consider a graph $\Gamma$  of complexity 1. 
Let $v_1$ denote the unique positive genus vertex.
Let $h=g(v_1)$.
Let $e$ be the degree of the unique edge incident to $v_1$.
Let $m=\{m_1, \ldots, m_l\}$ be the degrees of 
remaining edges of $\Gamma$.
The triple $(h,e,m)$ satisfies
$h\leq g$, $e\leq d$, and 
$m$ is a partition of $d-e$. 
The set of such triples
is in bijective correspondence to the set of
complexity 1 graphs:
$$(h,e,m)  \leftrightarrow \Gamma(h,e,m).$$
The contribution of $\Gamma(h,e,m)$ to $P(g,d)$ 
contains two Hodge integrals:  at the vertex $v_1$ and at the
vertex $v$ lying over $p_2$.
The Hodge integral at $v_1$ is $Q_h^e$ (up to signs).
The Hodge integral at $v$ is a $\lambda_g$ integral (see [FP2])
and may be integrated by the Virasoro constraints (\ref{lamg}).
A 
direct computation then
yields the contribution of $\Gamma$ to be:
\begin{equation}
\label{kdkdk}
\frac{(-1)^{d-g}} {d} \ Q^e_h\  \frac{e^{e+1}}{e!} \
\frac{(2h+d-l-1)!}{(2h-l)!} \
\frac{(-d)^l}{|\text{Aut}(m)|} \ \prod_{i=1}^l \frac{ m_i^{m_i-1}}{m_i!}
\end{equation}
$$
\cdot \  d^{2g-2h} \int_{\overline{M}_{g-h,1}} \psi_1^{2g-2h-2} \lambda_{g-h}.
$$
Here, $\text{Aut}(m)$ is the group which
permutes equal parts of $m$. The contribution vanishes unless
$2h\geq l$. Finally,
the integral
$\int_{\overline{M}_{0,1}} \psi_1^{-2} \lambda_0$ occurring
in (\ref{kdkdk}) in case $g=h$ is
defined to be 1. 

The integral $P(g,d)$ equals the sum of all graph
contributions (\ref{jjjj}--\ref{kdkdk}). 
As $P(g,d)=0$, we have found a relation among the Hodge integrals
including the $Q$ integrals.

\subsection{Proof of Proposition 1}
The Hodge relation found in Section 2.6 can be rewritten
using the following observations.

The Hodge 
integrals other than the $Q$ integrals appearing in 
(\ref{jjjj}--\ref{kdkdk}) are determined in [FP1]:
\begin{equation}
\label{bgcg}
\sum_{g\geq 0} d^{2g} t^{2g} \int_{\overline{M}_{g,1}} \psi_1^{2g-2} 
\lambda_{g}
 = \Big( \frac{dt/2}{\sin(dt/2)} \Big),
\end{equation}
\begin{equation}
\label{cccg}
\sum_{g\geq 1} d^{2g} t^{2g}
\int_{\overline{M}_{g,1}} \psi_1^{2g-1} \lambda_{g-1}
= 
\Big( \frac{dt/2}{\sin(dt/2)} \Big) 
\cdot \log\Big( \frac{dt/2}{\sin(dt/2)} \Big).
\end{equation}
Let $Part(a,b)$ denote the set of partitions of
$a$ of length $b$. 
The equality
$$f_{hde}=
\frac{e^{e+1}}{e!}
\sum_{l=0}^{2h} 
\frac{(2h+d-l-1)!}{(2h-l)!}
\sum_{m \in Part(d-e,l)}
\frac{(-d)^l}{|\text{Aut}(m)|} \ \prod_{i=1}^l \frac{ m_i^{m_i-1}}{m_i!}$$
follows directly from the definition (\ref{jack}).

Let $d\geq 1$ be fixed. The Hodge integral relations obtained
from the vanishing of $P(g,d)$ for all $g\geq 1$ may then be
expressed
as a series equality:
$$(\sum_{e=1}^d \sum_{g=1}^\infty Q_g^e f_{gde} t^{2g})
\ \cdot \ \Big( \frac{dt/2}{\sin(dt/2)} \Big) =
d^{d-1} 
\Big( \frac{dt/2}{\sin(dt/2)} \Big) 
\cdot \log\Big( \frac{dt/2}{\sin(dt/2)} \Big).$$
Proposition 1 follows from cancelling the invertible
series (\ref{bgcg}).

\section{\bf{Theorem 2}}
\subsection{Reduction}
The derivation of Theorem 2 from Proposition 1 requires
some knowledge of $\tau(x)$ and a significant amount of
binomial combinatorics.

Let $k$ be a fixed positive integer.
We start by summing the right side of (\ref{qwqw}) using Proposition 1:
\begin{equation}
\label{dfgh}
\sum_{j=1}^k { (-1)^{k-j}}
\frac{j^{k-1}}
{k} \binom {k}{j} 
 {\text {log}} \Big( \frac{jt/2}{{\text{sin}}(jt/2)} \Big)
\end{equation}
$$
=\sum_{g=1}^{\infty} t^{2g} \sum_{e=1}^k Q^e_g \sum_{j=e}^k
{ (-1)^{k-j}}
\frac{j^{k-j}}
{k} \binom {k}{j}f_{gje}. 
$$
A direct partial fraction expansion shows the equality:
$$
I(g,k) 
= \sum_{e=1}^k Q^e_g \  (-1)^{k-e} \frac{e^k}{k!} \binom {k}{e}.$$
Hence, Theorem 2 is a direct consequence of
(\ref{dfgh}) and the following Proposition.

\begin{pr}
Let $k\geq e$. Then,
$$\sum_{j=e}^k 
{ (-1)^{k-j}}
\frac{{j^{k-j}}}{k}
 \binom {k}{j}f_{gje} = \frac{(2g+k-1)!}{(2g)!} \cdot (-1)^{k-e}
  \frac{e^k}{k!} \binom {k} {e}.$$
\end{pr}

\subsection{Powers of $\tau$}
In order to prove Proposition 2, we will need a formula
for the coefficients of $\tau^l(x)$ appearing in the 
definition (\ref{jack}) of $f_{gje}$.
\begin{lm} Let $r,l \geq 0$,
\label{ccrr}
$$\frac{1}{l!} 
{\mathcal{C}}(x^r,\tau^l(x))=  \binom{r-1}{l-1} \frac{r^{r-l}}{r!} 
\ .$$
\end{lm}

\bpf This is a direct application of the Lagrange inversion formula
(see [dB], (2.2.4)).
Solving $x=z/f(z)$ with $f(z)=e^z$ gives 
$$z=\tau(x)=\sum_{r=1}^{\infty}
c_rx^r,$$ $$c_r=\frac{1}{r!}
[(d/dz)^{r-1}(f(z))^r]_{z=0}=r^{r-1}/(r!).$$
This is simply 
the well-known formula stated in Section 2.1. More generally, 
$$g(z)=g(0)+
\sum_{r=1}^{\infty}d_rx^r,$$  
$$d_r=\frac{1}{r!}[(d/dz)^{r-1}
\{g'(z)(f(z))^r\}]_{z=0}\,.$$
Applying this with $g(z)=z^l$ gives the result.
\epf

\subsection{Proof of Proposition 2}

Using definition (\ref{jack}), Lemma \ref{ccrr}, and simple
manipulations, we find Proposition 2 is equivalent to the
equation:
\begin{equation}
\label{hples}
\sum_{j=e+1}^k  \sum_{l=1}^{j-e}
\binom{2g+j-l-1}{j-1}\binom{k}{j} \binom{j-1}{e-1}
 \binom{j-e-1}{l-1}
j^{k-j+l} (e-j)^{j-e-l}
\end{equation}
$$= e^{k-e}\binom{k}{e} \Big( \binom{2g+k-1}{k-1} -
\binom{2g+e-1}{e-1} \Big).$$
To proceed, we may write the left and right sides of the
above equation canonically in terms of
the binomials
$$\binom{2g+e-1}{t+e-1}$$
for $0\leq t \leq k-e$ using the relations:
 $$\binom{2g+j-l-1}{j-1}= \sum_{t=l}^{j-e} \binom{j-e-l}{t-l} \binom
{2g+e-1}{t+e-1},$$
$$\binom{2g+k-1}{k-1} = \sum_{t=0}^{k-e} \binom{k-e}{t} 
\binom{2g+e-1}{t+e-1}.$$
Then it suffices to match the coefficients
\begin{equation}
\label{pwp}
\sum_{j=e+1}^k  \sum_{l=1}^{j-e} \binom{j-e-l}{t-l}\binom{k}{j}
\binom{j-1}{e-1}
 \binom{j-e-1}{l-1}
j^{k-j+l} (e-j)^{j-e-l}
\end{equation}
$$= e^{k-e} \binom{k}{e} \binom{k-e}{t}$$
for $1 \leq t \leq k-e$ (the matching at $t=0$ is
trivial).
Equation (\ref{pwp}) simplifies to:
$$
\sum_{j=e+1}^k  \sum_{l=1}^{j-e} \binom{k-e-t}{j-e-t}
\binom{t-1}{l-1}
j^{k-j-1+l} (e-j)^{j-e-1-l}
= -\frac{ e^{k-e-1} }{t}.
$$
Summing over $l$ yields:
$$
\sum_{j=e+t}^k  \binom{k-e-t}{j-e-t}
j^{k-j} (e-j)^{j-e-t-1}
= -\frac{ e^{k-e-t} }{t}.
$$
Substitute $z=k-e$, $s=j-e-t$. Then, we must prove
\begin{equation}
\label{weww}
\sum_{s=0}^{z-t}  \binom{z-t}{s}
(e+s+t)^{z-t-s} (-s-t)^{s-1}
= -\frac{ e^{z-t} }{t},
\end{equation}
for all $1 \leq t \leq z$.
If the left side of (\ref{weww}) is viewed as a polynomial
in $e$, the coefficient of $e^{z-t}$ clearly matches the right side.
Hence, it suffices to show the coefficient
of $e^q$ vanishes for $0\leq q <z-t$:
$$\sum_{s=0}^{z-t-q} \binom{z-t}{s} \binom{z-t-s}{q}(s+t)^{z-t-s-q}  
 (-s-t)^{s-1}=0.$$
This is equivalent to:
$$\sum_{s=0}^{z-t-q} \binom{z-t-q}{s}(s+t)^{z-t-s-q}  
 (-s-t)^{s-1}=0.$$
Substituting $n=z-t-q$ and simplifying, we must prove:
\begin{equation}
\label{kfgg}
\sum_{s=0}^{n} (-1)^s \binom{n}{s}(s+t)^{n-1}  
=0,
\end{equation}
for all $n>0$. 
Finally, the proof of Proposition 2 (and therefore
Theorem 2) is completed by observing
(\ref{kfgg}) follows from the well-known relation:
$$\sum_{s=0}^n (-1)^s \binom{n}{s} s^\gamma = 0$$
for all $0\leq \gamma \leq n-1$.
\epf

\section{\bf{Hyperelliptic Hodge integrals}}
In this section we compute for all $g$ the $M_g$-evaluation
of the class of the hyperelliptic locus $H_g$. As explained
in the Introduction, this provides an alternative proof
of Theorem 1 in the case $k=2$ and its Corollary 2.

As in Section 1, the starting point is the identity
\begin{equation} \label{start}
\lambda_{g}\lambda_{g-1}=(-1)^{g-1}(2g-1)!\,\ch_{2g-1}(\hodge).
\end{equation}
Mumford's calculation of the Chern character of the Hodge
bundle [Mu] gives then an expression for
$\lambda_{g}\lambda_{g-1}$ in terms of $\kappa$ and $\psi$ classes.
This expression lends itself very well for a direct evaluation
on the hyperelliptic locus: in the usual model
of hyperelliptic curves as double covers of rational curves, all
relevant classes are pullbacks from the moduli of rational curves, where
evaluation is straightforward. In the process one finds simple
expressions (in the rational model) for all components of the
restriction of $\ch(\hodge)$ to the hyperelliptic locus.
This generalizes the formula of Cornalba and Harris [CH]
for $\lambda_1$ on $\overline{H}_g$.
It seems plausible that these expressions will allow the evaluation
of other hyperelliptic Hodge integrals. 

We may view $\overline{M}_{0,2g+2}$ as the coarse
moduli space of stable 
hyperelliptic curves of genus $g$
with an ordering of the Weierstrass points
(see [HM] 6C or [FP1] \S3.2). 
The universal hyperelliptic curve is then the (stack) double cover
of $\overline{M}_{0,2g+3}$ branched over $B$, the
disjoint union of the $2g+2$ sections:
\begin{equation*}
\begin{CD}
 \UM
@>{f}>>  \overline{M}_{0,2g+3}\, \\
@V{\varpi}VV   @V{\pi}VV \\
 \overline{H}_g^{{\rm ord}} @= \,\overline{M}_{0,2g+2}\,.
\end{CD}
\end{equation*}
We have $\psi_1=f^*(\psi_{2g+3}-B/2)$.
Writing $h_i$ for the genus $g$ class $\kappa_i$ viewed on
$\overline{M}_{0,2g+2}\,$, we obtain:
\begin{eqnarray*}
h_i&=&\varpi_*\psi_1^{i+1}=\varpi_*(f^*(\psi_{2g+3}-B/2))^{i+1}
=\pi_*f_*f^*((\psi_{2g+3}-B/2)^{i+1}) \\
&=& 2\pi_*((\psi_{2g+3}-B/2)^{i+1})
= 2\pi_*(\psi_{2g+3}^{i+1}+(-B/2)^{i+1})        \\
&=& 2\kappa_i + 2\sum_{j=1}^{2g+2}(-\tfrac12)^{i+1}(-\psi_j)^i 
= 2\kappa_i - 2^{-i}\sum_{j=1}^{2g+2}\psi_j^i\,.
\end{eqnarray*}
(Here the genus 0 class $\kappa_i$ in the last line is the
generalization to $\overline{M}_{g,n}$ by Arbarello-Cornalba [AC]
of Mumford's class for $\overline{M}_g$.)
Writing $\chi_i=\ch_i(\hodge)$, we have computed the first term
in Mumford's formula 
$$
\frac{(2k)!}{B_{2k}}\,\chi_{2k-1}=\kappa_{2k-1}
+\frac12\sum_{h=0}^{g-1}i_{h,*}
\frac{\psi_1^{2k-1}+\psi_2^{2k-1}}{\psi_1+\psi_2}
$$
in the rational model, and it remains to evaluate the boundary
terms. (Recall that $\chi_{2k}=0$ for positive $k$.)

Boundary divisors of $\overline{M}_{0,2g+2}$ come in two types:
odd boundary divisors, with an underlying partition of $2g+2$
in two odd numbers ($\ge3$), and even boundary divisors. As
described in [CH] and [HM], the hyperelliptic curves
corresponding to an odd boundary divisor generically have one
disconnecting node and four automorphisms, while those
corresponding to an even boundary divisor generically have two
non-disconnecting nodes and two automorphisms.

As a result, Mumford's formula in codimension one reads on the
rational model as follows:
$$ 12\chi_1=2\kappa_1-\tfrac12\psi+\tfrac12\delta_{{\rm odd}}
+2\delta_{{\rm even}}
$$
with evident notations. Since $\kappa_1=\psi-\delta$ in genus 0,
this simplifies to
$$ 8\chi_1=\psi-\delta_{{\rm odd}}=\kappa_1+\delta_{{\rm even}}\,.
$$
The higher codimension case is very similar.
The terms with $1\le h\le g-1$ in Mumford's formula correspond
to the odd boundary divisors. In the rational model they appear
with an extra factor $\frac12$. Now $\psi_1=f_h^*(\psi_{2h+3}-B/2)$;
since this is here a cotangent line at a Weierstrass point,
we must evaluate $\psi_{2h+3}-B/2$ on a Weierstrass point divisor
in $\overline{M}_{0,2h+3}$. It is easy to check that the result,
as a class on a boundary divisor of $\overline{M}_{0,2g+2}$
with underlying partition $[2h+1,2(g-h)+1]$, is $\frac12\psi_*$,
where $\psi_*$
is the cotangent line in the node to the branch with $2h+1$
marked points. Analogously, for $\psi_2$ and genus $g-h$,
we find $\frac12\psi_\bullet$, where $\psi_\bullet$ is the
cotangent line in the node to the other branch.
Therefore the odd boundary contribution 
to $\frac{(2k)!}{B_{2k}}\,\chi_{2k-1}$ equals
$$
\frac12\sum_{{\rm odd}\,D}\frac{(\frac12\psi_*)^{2k-1}
+(\frac12\psi_\bullet)^{2k-1}}
{\frac12\psi_*+\frac12\psi_\bullet}\bigg|_D
=\frac1{2^{2k-1}}\sum_{{\rm odd}\,D}\frac{\psi_*^{2k-1}
+\psi_\bullet^{2k-1}}
{\psi_*+\psi_\bullet}\bigg|_D\,.
$$
The $h=0$ term in Mumford's formula breaks up in terms
corresponding to the even boundary divisors; each of these
appears with an extra factor 2. To identify the classes
$\psi_1$ and $\psi_2\,$, we need to construct the family
of hyperelliptic curves corresponding to an even boundary
divisor with underlying partition $[2h+2,2k+2]$
(hence $h+k=g-1$). The base of the family is
$\UM_h\times\UM_k\,$. The idea is to glue 
$\UM_h\times_{\overline{H}_h}\UM_h$ and
$\UM_k\times_{\overline{H}_k}\UM_k$
along two sections on either side, the diagonal $\Delta$
and its image $\Delta'=\{(p,p')\}$ under the hyperelliptic
involution on the second factor. However, $\Delta$ and
$\Delta'$ intersect along $\Delta(W)$, where $W$ is the
Weierstrass divisor in $\UM$. Therefore
$\UM\times_{\overline{H}}\UM$
must be blown up along $\Delta(W)$, on either side.
The relative canonical divisor induced on the second factor
after the blow-up can be identified with the class
$\psi_1+W$ on the second factor before blowing up.
Therefore the classes $\psi_1$ and $\psi_2$ in Mumford's
formula correspond on the rational model to
$f_h^*(\psi_{2h+3})$ and $f_k^*(\psi_{2k+3})$ respectively,
and the even boundary contribution to 
$\frac{(2k)!}{B_{2k}}\,\chi_{2k-1}$ equals simply
$$2\sum_{{\rm even}\,D}\frac{\psi_*^{2k-1}
+\psi_\bullet^{2k-1}}
{\psi_*+\psi_\bullet}\bigg|_D\,.
$$
We have proven:

\begin{pr} \label{hyphodge}
In the coarse rational model $\overline{M}_{0,2g+2}=
\overline{H}_g^{{\rm ord}}$,
the Chern character of the genus $g$ Hodge bundle equals
\begin{eqnarray*}
\lefteqn{
{\rm ch}(\hodge)=g+\sum_{k=1}^g\frac{B_{2k}}{(2k)!}
\left[2\kappa_{2k-1}-\frac1{2^{2k-1}}\sum_{j=1}^{2g+2}
\psi_j^{2k-1} \right.   } \cr
&&\qquad\qquad \left. \mbox{}+\frac1{2^{2k-1}}
\sum_{{\rm odd}\,D}\frac{\psi_*^{2k-1}
+\psi_\bullet^{2k-1}}
{\psi_*+\psi_\bullet}\bigg|_D
+2\sum_{{\rm even}\,D}\frac{\psi_*^{2k-1}
+\psi_\bullet^{2k-1}}
{\psi_*+\psi_\bullet}\bigg|_D\right].
\end{eqnarray*}
\end{pr}
\noindent(The vanishing of $\ch(\hodge)$ in degrees $\ge2g$
--- here trivial --- holds on $\overline{M}_g$ as well,
see e.g.~[FP1].)

In fact, these formulas can be simplified, just as in
codimension 1:
\begin{eqnarray*}
\frac{(2k)!}{B_{2k}}\,\chi_{2k-1} & = &
\frac{2^{2k}-1}{2^{2k-1}}
\left(\sum_{j=1}^{2g+2}
\psi_j^{2k-1}-\sum_{{\rm odd}\,D}\frac{\psi_*^{2k-1}
+\psi_\bullet^{2k-1}}
{\psi_*+\psi_\bullet}\bigg|_D \right) \\
& = &
\frac{2^{2k}-1}{2^{2k-1}}
\left( \kappa_{2k-1} +
\sum_{{\rm even}\,D}\frac{\psi_*^{2k-1}
+\psi_\bullet^{2k-1}}
{\psi_*+\psi_\bullet}\bigg|_D \right).
\end{eqnarray*}
This follows from the identity
$$\kappa_{2k-1}
= \sum_{j=1}^n \psi_j^{2k-1}
- \frac{\psi_*^{2k-1}
+\psi_\bullet^{2k-1}}
{\psi_*+\psi_\bullet}\bigg|_\delta
$$
on $\mbar{0,n}\,$, a consequence of Proposition 1 in [FP1].

\bigskip
\noindent{\bf Corollary.}
\label{hypcor}
{\it
On $\overline{H}_g^{{\rm ord}}$,
$${\rm ch}_{2g-1}(\hodge)
= \frac{B_{2g}}{(2g)!}\,(2^{2g+1}-2).
$$
Hence on the stack $\overline{H}_g\,$,
}
$$
\lmb{g}\lmb{g-1} =
\frac{(2^{2g}-1)|B_{2g}|}{(2g+2)!\,2g}\,.
$$

\bpf By the above
$$
\frac{(2g)!}{B_{2g}}\,
\chi_{2g-1}= \frac{2^{2g}-1}{2^{2g-1}}\left(1+\frac12
\sum_{h=1}^g \binom{2g+2}{2h} \right)
= \frac{2^{2g}-1}{2^{2g-1}}\,2^{2g}=2^{2g+1}-2,
$$
whence the first formula. The second formula follows by 
using (\ref{start}) and dividing
by $2\cdot(2g+2)!\,$. The factor of 2 is required to account for
the hyperelliptic automorphism groups in
the stack $\overline{H}_g\,$.
\epf

\section{\bf{Theorem 2 revisited}}
\subsection{Reformulation}
In this section we present a reformulation of Theorem 2
that reduces all known (and several conjectured)
non-vanishing results to combinatorial identities.
For $g\ge1$, consider the polynomial $P_g(k)$ in $k$ of degree $g-1$ 
(with zero constant term for $g\ge2$) defined by:
$$
\frac{|B_{2g}|}{2g}\,
P_g(k)=
\sum_{i=0}^{g-1}(-1)^ik^{g-1-i}
\int_{\overline{M}_{g,1}} \psi_1^{g-1-i}\lmb{i}\lmb{g}\lmb{g-1} \,.
$$
Note that the right-hand side equals $Q_g^k$ as in (\ref{mqz})
for positive integers $k$.

\vspace{+10pt}
\noindent{\bf Theorem 3.}
{\em
For positive integers $g$,  $k$,}
$$
P_g(k)=
\sum_{l=1}^k
\frac{(k-1)!}{(k-l)!}\,\frac1{k^l}
\sum_{m=1}^l(-1)^{l-m}
\binom{l}{m}
\frac{m^{2g+l-1}}{(2g+l-1)!}\,.
$$
\vspace{+10pt}

\bpf
This follows directly from Theorem 2.
By expanding the logarithmic series as in [FP1], Lemma 3,
one obtains
$$
I(g,k)=\frac{(k-1)!}{(2g+k-1)!}\frac{|B_{2g}|}{2g}
\sum_{j=1}^k\frac{(-1)^{k-j}}{(k-j)!}\frac{j^{k-1}}{j!}j^{2g}\,.
$$
Since
$$
\frac1{\prod_{i=1}^k(1-i\psi_1)}=\sum_{n=0}^{\infty}\psi_1^n
\frac{(-1)^k}{k!}\sum_{j=1}^k(-1)^jj^{k+n}\binom{k}j
$$
we also have
$$
I(g,k)=\int_{\overline{M}_{g,1}} \lmb{g}\lmb{g-1}c(\hodge^*)
\sum_{n=0}^{\infty}\psi_1^n\frac1{k!}\sum_{j=1}^k
(-1)^{k-j}j^{k-1}\binom{k}{j}j^{n+1}\,.
$$
Now observe that the resulting identity can be written
as $BA=DBV$, where $A$ is the infinite vector with entries
$$A(j)=\int_{\overline{M}_{g,1}} \lmb{g}\lmb{g-1}c(\hodge^*)
\sum_{n=0}^{g-1}j^{n+1}\psi_1^n$$ (for a fixed $g$),
$B$ is the infinite lower-triangular matrix with entries
$$B(i,j)=(-1)^{i+j}j^{i-1}\binom{i}j$$ for $1\le j\le i$,
$D$ is the infinite diagonal matrix with entries
$$D(k,k)=\frac{(k-1)!}{(2g+k-1)!}\frac{|B_{2g}|}{2g},$$
and $V$ is the infinite vector with entries $V(j)=j^{2g}$.

One easily shows that the inverse of $B$ has entries
$B^{-1}(i,j)=\binom{i-1}{j-1}i^{1-j}$ for $1\le j\le i$.
The Theorem follows by writing
out $A=B^{-1}DBV$ and using $\frac{|B_{2g}|}{2g}\,P_g(k)=A(k)/k$.
\epf

The connection to the Stirling number formula in Section \ref{sss}
is obtained from the equation:
$${\mathfrak{S}}^{(l)}_{2g-1+l}= \frac{1}{l!} \sum_{m=1}^l(-1)^{l-m}
\binom{l}{m}
{m^{2g+l-1}}.$$

\subsection{Non-vanishing results}
We present here the reformulations of four non-vanishing results.
All four are proved by D.~Zagier in the Appendix from Theorem 3.
Equivalently, these are 
identities in the socle of the tautological ring $R^*(M_g)$.
First, the leading coefficient in $P_g(k)$ is:
\begin{equation}
\label{eg-1}
 {\mathcal{C}}(k^{g-1}, P_g(k)) = \frac1{2^{2g-1}(2g-1)!!}\,.
\end{equation}
Equation (\ref{eg-1}) is equivalent to Theorem 1 (providing
an alternate proof which avoids the KdV equations).
The next highest coefficient is:
\begin{equation}
\label{eg-2}
{\mathcal{C}}(k^{g-2}, P_g(k)) =
\frac{-g(g-2)}{3^22^{2g-1}(2g-1)!!}\,,
\end{equation}
in agreement with the prediction for $\kp{g-3}\lmb1$ in [F3].
Zagier has found generalizations of these combinatorial
formulas for the coefficient of $k^{g-1-i}$ in $P_g(k)$  (for
fixed codegree $i$).

Similarly, Bernoulli number formulas are found 
in the Appendix for the coefficient
of $k^i$ in $P_g(k)$
for fixed degree $i$.
The coefficient of the linear term in $P_g(k)$ is:
\begin{equation}
\label{eg1}
{\mathcal{C}}(k^{1}, P_g(k)) =
\frac{B_{2g-2}}{2\cdot(2g-2)!}\,,
\end{equation}
in agreement with (\ref{lglglg}) previously calculated in [FP1].
The quadratic coefficient in $P_g(k)$ is:
\begin{equation}
\label{eg2}
{\mathcal{C}}(k^{2}, P_g(k))=
\frac{-g\,B_{2g-2}}{2\cdot(2g-2)!}\,.
\end{equation}
Equation (\ref{eg2}) determines the evaluation of $\kp1\lmb{g-3}$
for $g\ge3$ --- it implies Conjecture 2 in [F4].

\vspace{+10 pt}
\noindent Department of Mathematics \hfill Institutionen f\"or Matematik \\
\noindent Oklahoma State University \hfill Kungliga Tekniska H\"ogskolan \\
\noindent Stillwater, OK 74078 \hfill 100 44 Stockholm, Sweden \\
\noindent cffaber@math.okstate.edu \hfill carel@math.kth.se


\vspace{+10 pt}
\noindent
Department of Mathematics \\
\noindent California Institute of Technology \\
\noindent Pasadena, CA 91125 \\
\noindent rahulp@cco.caltech.edu

\begin{thebibliography}{[COGP]}

\bibitem[AC]{AC} E. Arbarello and M. Cornalba, {\em
Combinatorial and algebro-geometric cohomology
classes on the moduli spaces of curves}, 
J. Algebraic Geom. 5 (1996), 705--749.

\bibitem[dB]{dB} N.G. de Bruijn, {\em
Asymptotic Methods in Analysis},
North-Holland/P.~Noordhoff, Amsterdam/Groningen 1958.

\bibitem[CH]{CH} M. Cornalba and J. Harris, {\em
Divisor classes associated to families of stable varieties,
with applications to the moduli space of curves}, 
Ann. Sci. \'Ecole Norm. Sup. (4) 21 (1988), no. 3, 455--475.

\bibitem[Dij]{Dij} R. Dijkgraaf, {\em
Some facts about tautological classes},
private communication, November 1993.

\bibitem[EHX]{ehx} T. Eguchi, K. Hori, and C.-S. Xiong, {\em
Quantum cohomology and Virasoro algebra}, Phys. Lett. {\bf B402} (1997),
71--80.

\bibitem[F1] {f1} C. Faber, {\em Chow rings of moduli spaces of curves. 
I. The Chow ring of $\overline{M}_3$}, 
Ann. of Math. (2) 132 (1990), no. 2,
331--419. 

\bibitem[F2] {f2} C. Faber, {\em Chow rings of moduli spaces 
of curves. II. Some results on the 
Chow ring of $\overline{M}_4$}, Ann. of Math. (2) 132
(1990), no. 3, 421--449. 

\bibitem[F3] {f3} C. Faber, {\em A conjectural description of the 
tautological ring of the moduli space of curves},
in {\sl Moduli of Curves and Abelian Varieties
(The Dutch Intercity Seminar on Moduli)\/} (C.~Faber and
E.~Looijenga, eds.),
109--129, Aspects of Mathematics E~33, Vieweg, Wiesbaden 1999.

\bibitem[F4]{f4} C. Faber, {\em
Algorithms for computing intersection numbers on moduli spaces of curves,
with an application to the class of the locus of Jacobians},
in New Trends in Algebraic Geometry
(K.~Hulek, F.~Catanese, C.~Peters and M.~Reid, eds.),
93--109, Cambridge University Press, 1999.

\bibitem[FL]{fl} C. Faber and E. Looijenga, {\em
Remarks on moduli of curves},
in {\sl Moduli of Curves and Abelian Varieties
(The Dutch Intercity Seminar on Moduli)\/} (C.~Faber and
E.~Looijenga, eds.),
23--45, Aspects of Mathematics E~33, Vieweg, Wiesbaden 1999.

\bibitem[FP1] {fp1} C. Faber and R. Pandharipande, {\em Hodge
integrals and Gromov-Witten theory}, 
Invent. Math. {\bf 139} (2000) 1, 173--199.

\bibitem[FP2] {fp2} C. Faber and R. Pandharipande, {\em
Hodge integrals, partition matrices, and the $\lambda_g$ conjecture},
preprint 1999, {\tt math.AG/9908052}.

\bibitem[FP3] {fp3} C. Faber and R. Pandharipande, in preparation.

\bibitem[FanP] {fanp} B. Fantechi and R. Pandharipande, {\em
Stable maps and branch divisors}, preprint 1999,
{\tt math.AG/9905104}.

\bibitem[FuP]{fup} W. Fulton and R. Pandharipande, 
{\em Notes on stable maps and quantum cohomology}, in
Proceedings of Symposia in Pure  Mathematics: Algebraic Geometry 
Santa Cruz 1995, J. Koll\'ar, R. Lazarsfeld, D. Morrison, eds.,
Volume 62, Part 2, p.45--96.

\bibitem[Ge]{Ge} E. Getzler, {\em
The Virasoro conjecture for Gromov-Witten invariants},
to appear in Algebraic Geometry --- Hirzebruch 70
(P. Pragacz et al., eds.), AMS, Contemporary Mathematics.

\bibitem[GeP]{gep} E. Getzler and R. Pandharipande, {\em Virasoro
constraints and the Chern classes of the Hodge bundle},
Nucl. Phys. {\bf B530} (1998),  701--714. 

\bibitem[GrP]{gp} T. Graber and R. Pandharipande, {\em Localization
of virtual classes}, Invent. Math. {\bf 135} (1999), 487--518.

\bibitem[HL]{hl} R. Hain and E. Looijenga, {\em Mapping class groups and
moduli spaces of curves}, in
Proceedings of Symposia in Pure  Mathematics: Algebraic Geometry 
Santa Cruz 1995, J. Koll\'ar, R. Lazarsfeld, D. Morrison, eds.,
Volume 62, Part 2, p.97--142.

\bibitem[HM]{hm} J. Harris and I. Morrison, {\em
Moduli of Curves}, GTM 187, Springer, 1998.

\bibitem[I]{i} E. Izadi, {\em
The Chow ring of the moduli space of curves of genus $5$},
in {\em The moduli space of curves},
(R. Dijkgraaf, C. Faber, and G. van der Geer, eds.), Birkh\"auser,
1995, 401--417.

\bibitem[K1] {k} M. Kontsevich, {\em Intersection theory on the moduli
space of curves and the matrix Airy function}, Comm. Math. Phys.
{\bf 147} (1992), 1--23.

\bibitem[K2]{k2}  M. Kontsevich, {\em Enumeration of rational curves via
torus actions}, in {\em The moduli space of curves}
(R. Dijkgraaf, C. Faber, and G. van der Geer, eds.), Birkh\"auser,
1995, 335--368.

\bibitem[L]{l} E. Looijenga, {\em On the tautological
ring of $M_g$}, Invent. Math. ${\mathbf{121}}$ (1995), 411--419.

\bibitem[Ma]{m} Yu. Manin, {\em Generating functions in
algebraic geometry and sums over trees},
in {\em The moduli space of curves},
(R. Dijkgraaf, C. Faber, and G. van der Geer, eds.), Birkh\"auser,
1995, 401--417.

\bibitem[Mu] {mu} D. Mumford, {\em Towards an enumerative geometry of
the moduli space of curves}, in {\em Arithmetic and Geometry}
(M. Artin and J. Tate, eds.), Part II, Birkh\"auser, 1983, 271--328.

\bibitem[P1]{pp} R. Pandharipande, {\em The Chow ring of the
non-linear Grassmannian}, J. Alg. Geom. ${\mathbf{7}}$ (1998), 123-140.

\bibitem[P2]{p} R. Pandharipande, {\em Hodge integrals and
degenerate contributions}, Comm. Math. Phys. (to appear).


\bibitem[W]{w2} E. Witten, {\em Two dimensional gravity and intersection
theory on moduli space}, Surveys in Diff. Geom. {\bf 1} (1991), 243--310.

\end{thebibliography}
\end{document}